\documentclass[11pt]{article}
\usepackage{hyperref}
\usepackage{amsfonts, amsmath, amssymb, amsthm, tikz, tkz-graph, enumitem, microtype}

\definecolor{darkgray}{RGB}{64,64,64}
\definecolor{litegray}{RGB}{192,192,192}
\definecolor{red}{rgb}{.859375,.265625,.21484375}
\definecolor{blue}{rgb}{.26171875,.51953125,.95703125}
\tikzstyle{vertex}=[circle, draw, fill=litegray, inner sep=0pt, minimum size=4pt]
\tikzstyle{bvertex}=[circle, draw, fill=litegray, inner sep=0pt, minimum size=6pt]
\tikzstyle{block}=[draw, rectangle, minimum height=3cm, minimum width=0.5cm, text centered, draw=darkgray, font=\small, inner sep=0]
\tikzstyle{block_h}=[draw, rectangle, minimum height=0.5cm, minimum width=3cm, text centered, draw=darkgray, font=\small, inner sep=0]

\usepackage{mathtools}

\DeclarePairedDelimiter{\ceil}{\lceil}{\rceil}

\title{Cooperative colorings of trees and of bipartite graphs}
\author{
Ron Aharoni\thanks{Department of Mathematics, Technion -- Israel Institute of Technology, Technion City, Haifa 3200003, Israel. Email: {\tt ra@tx.technion.ac.il}. Supported in part by the United States--Israel Binational Science Foundation (BSF) grant no.\ 2006099, the Israel Science Foundation (ISF) grant no.\ 2023464 and the Discount Bank Chair at Technion. This paper is part of a project that has received funding from the European Union's Horizon 2020 research and innovation programme, under the Marie Sk{\l{}}odowska-Curie grant agreement no.\ 823748.}
\and Eli Berger\thanks{Department of Mathematics, University of Haifa, Mt. Carmel, Haifa 3498838, Israel. Email: {\tt berger@math.haifa.ac.il}. Supported in part by BSF grant no.\ 2006099 and ISF grant no.\ 2023464.}
\and Maria Chudnovsky\thanks{Mathematics Department, Princeton University, Princeton, NJ 08544, USA. Email: {\tt mchudnov@math.princeton.edu}. Supported in part by BSF grant no.\ 2006099, NSF grant DMS-1550991 and US Army Research Office Grant W911NF-16-1-0404.}
\and Fr\'ed\'eric Havet\thanks{CNRS, Universit\'e C\^ote d'Azur, I3S, and INRIA, Sophia-Antipolis Cedex 06902, France. Email: {\tt frederic.havet@inria.fr}.}
\and Zilin Jiang\thanks{Department of Mathematics, Massachusetts Institute of Technology, Cambridge, MA 02139, USA. Email: {\tt zilinj@mit.edu}. The work was done when Z. Jiang was a postdoctoral fellow at Technion -- Israel Institute of Technology, and was supported in part by ISF grant nos.\ 409/16, 936/16.}}
\date{}

\newtheorem{theorem}{Theorem}

\theoremstyle{definition}
\newtheorem{definition}{Definition}

\theoremstyle{remark}
\newtheorem{remark}{Remark}

\newcommand{\cb}{\mathcal{B}}
\newcommand{\cc}{\mathcal{C}}
\newcommand{\cf}{\mathcal{F}}
\newcommand{\cg}{\mathcal{G}}

\newcommand{\ci}{\mathcal{I}}
\newcommand{\cm}{\mathcal{M}}
\newcommand{\cp}{\mathcal{P}}
\newcommand{\cs}{\mathcal{S}}
\newcommand{\ct}{\mathcal{T}}

\newcommand{\N}{\mathbb{N}}

\newcommand{\from}{\colon}
\newcommand{\sset}[1]{\left\{#1\right\}}
\newcommand{\dset}[2]{\left\{#1:#2\right\}}
\newcommand{\abs}[1]{\left\lvert#1\right\rvert}
\newcommand{\pr}[1]{\operatorname{Pr}\left(#1\right)}
\newcommand{\cpr}[2]{\pr{#1 \mid #2}}

\begin{document}

\maketitle

\begin{abstract}
Given a system $(G_1, \ldots ,G_m)$ of graphs on the same vertex set $V$, a cooperative coloring is a choice of vertex sets $I_1, \ldots ,I_m$, such that $I_j$ is independent in $G_j$ and $\bigcup_{j=1}^{m}I_j = V$. For a class $\mathcal{G}$ of graphs, let $m_{\mathcal{G}}(d)$ be the minimal $m$ such that every $m$ graphs from $\mathcal{G}$ with maximum degree $d$ have a cooperative coloring. We prove that $\Omega(\log\log d) \le m_\mathcal{T}(d) \le O(\log d)$ and $\Omega(\log d)\le m_\mathcal{B}(d) \le O(d/\log d)$, where $\mathcal{T}$ is the class of trees and $\mathcal{B}$ is the class of bipartite graphs.
\end{abstract}

\section{Introduction}

A set of vertices in a graph is called \emph{independent} if no two vertices in it form an edge. A \emph{coloring} of a graph $G$ is a covering of $V(G)$ by independent sets. Given a system $(G_1, \ldots ,G_m)$ of graphs on the same vertex set $V$, a \emph{cooperative coloring} is a choice of vertex sets $\dset{I_j \subseteq V}{j \in [m]}$ such that $I_j$ is independent in $G_j$ and $\bigcup_{j=1}^{m}I_j =V$. If all $G_j$'s are the same graph $G$, then a cooperative coloring is just a proper vertex coloring of $G$ by $m$ independent sets.

A basic fact about vertex coloring is that every graph $G$ of maximum degree $d$ is $(d+1)$-colorable. It is therefore natural to ask whether $d+1$ graphs, each of maximum degree $d$, always have a cooperative coloring. This was shown to be false:

\begin{theorem}[Theorem 5.1 of Aharoni, Holzman, Howard and Spr\"ussel~\cite{MR3359930}]
  For every $d \ge 2$, there exist $d+1$ graphs of maximum degree $d$ that do not have a cooperative coloring.
\end{theorem}

% A cooperative coloring problem can be formulated in terms of independent transversals (this will be elaborated on in Section~\ref{sec:covers}).
Using the fundamental result on independent transversals of Haxell~\cite[Theorem~2]{MR1860440}, it can be shown that $2d$ graphs of maximum degree $d$ always have a cooperative coloring. Let $m(d)$ be the minimal $m$ such that every $m$ graphs of maximum degree $d$ have a cooperative coloring. By the above, $m(1) = 2$ and
\begin{equation}\label{basic}
  d+2 \le m(d) \le 2d, \text{ for every }d\ge 2.
\end{equation}
The theorem of Loh and Sudakov~\cite[Theorem 4.1]{MR2354708} on independent transversals in locally sparse graphs implies that $m(d) = d + o(d)$. Neither the lower bound nor the upper bound in \eqref{basic} has been improved for general $d$; even $m(3)$ is not known. However, restricting the graphs to specific classes, better upper bounds can be obtained.

\begin{definition}
  For a class $\cg$ of graphs, denote by $m_\cg(d)$ the minimal $m$ such that every $m$ graphs belonging to $\cg$, each of maximum degree at most $d$, have a cooperative coloring.
\end{definition}

For example, the following was proved:
\begin{theorem}[Corollary 3.3 of Aharoni et al.~\cite{MR2345810} and Theorem 6.6 of Aharoni et al.~\cite{MR3359930}]
Let $\cc$ be the class of chordal graphs and let $\cp$ be the class of paths. Then $m_{\cc}(d)=d+1$ for all $d$, and $m_{\cp}(2) = 3$.
\end{theorem}

In this paper, we prove some bounds on $m_{\cg}(d)$ for two more classes:

\begin{theorem}\label{main}
  Let $\ct$ be the class of trees, and let $\cb$ be the class of bipartite graphs. Then for $d \ge 2$,
  \begin{align*}
    \log_2\log_2d \le & m_{\ct}(d) \le \left(1+o(1)\right)\log_{4/3}d, \\
    \log_2 d\le & m_\cb(d) \le (1+o(1)) \frac{2d}{\ln d}.
  \end{align*}
\end{theorem}

\begin{remark}
  Let $\cf$ be the class of forests. It is evident that $m_\cf(d) \ge m_\ct(d)$ as $\cf \supset \ct$. Conversely, when $d \ge 2$, given $m = m_\ct(d)$ forests $F_1, \dots, F_m$ of maximum degree $d$, we can add edges to $F_i$ to obtain a tree $F_i'$ of maximum degree $d$, and the cooperative coloring for $F_1', \dots, F_m'$ is also a cooperative coloring for $F_1, \dots, F_m$. Therefore $m_\cf(d) = m_\ct(d)$ for $d \ge 2$.
\end{remark}

The notions of cooperative coloring and of list coloring have a common generalization: given a system $(G_1, \dots, G_m)$ of graphs with vertex sets $V_1, \dots, V_m$ (which are not neccessarily the same vertex set), a \emph{cooperative list coloring} is then a choice of independent sets in $G_i$ whose union equals $V := V_1 \cup \dots \cup V_m$. The notion of cooperative coloring is obtained by taking $V_i=V$, and list colorings are formed when $G_i$ is an induced subgraph of the same graph $G$ for all $i$. The upper bounds in Theorem~\ref{main} generalize to cooperative list colorings. For example, our proof of Theorem~\ref{main} for bipartite graphs readily gives the following result.

\begin{theorem} \label{thm:clc}
  For every system $(G_1, \dots, G_m)$ of bipartite graphs with maximum degree $d$ with vertex sets $V_1, \dots, V_m$, there is a cooperative list coloring if for every $v \in V_1 \cup \dots \cup V_m$, the number of its occurrences in $V_1, \dots, V_m$, that is $\abs{\dset{i \in [m]}{v \in V_i}}$, is at least $(1+o(1))\frac{2d}{\ln d}$.
\end{theorem}

A conjecture of Alon and Krivelevich~\cite[Conjecture~5.1]{MR1774970} states that the choice number of any bipartite graph with maximum degree $d$ is at most $O(\log d)$ (see \cite{MR2570629} for a result in this direction). This conjecture would follow if the term $(1+o(1))\frac{2d}{\ln d}$ in Theorem~\ref{thm:clc} was strengthened to $\Omega(\log d)$.

The rest of the paper is organized as follows. In Section~\ref{sec:trees} and Section~\ref{sec:bipartite}, we prove Theorem~\ref{main} for trees and bipartite graphs respectively. In Section~\ref{sec:covers} we discuss a further generalization of cooperative colorings.

\section{Trees} \label{sec:trees}

\begin{proof}[Proof of the lower bound on $m_{\ct}(d)$]
  Note that the system $\ct_2$, consisting of two paths in Figure~\ref{paths2} (one in thin red, the other in bold blue), does not have a cooperative coloring.

  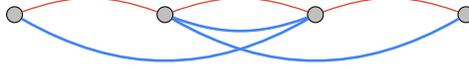
\begin{figure}[t] \centering
    \begin{tikzpicture}[scale=1]
      \node[bvertex] (a) at (0,0) {};
      \node[bvertex] (b) at (2,0) {};
      \node[bvertex] (c) at (4,0) {};
      \node[bvertex] (d) at (6,0) {};

      \path[-,color=red,line width=0.5pt]
        (a) edge[bend left=20] (b)
        (b) edge[bend left=20] (c)
        (c) edge[bend left=20] (d);

      \path[-,color=blue,line width=1pt]
        (a) edge[bend left=-30] (c)
        (b) edge[bend left=-20] (c)
        (b) edge[bend left=-30] (d);
    \end{tikzpicture}
    \caption{Construction of two paths without a cooperative coloring.} \label{paths2}
  \end{figure}

  Suppose now that $\cs=(F_1, F_2, \dots, F_m)$ is a system of forests on a vertex set $V$, not having a cooperative coloring. We shall construct a system $Q(\cs)$ of $m+1$ new forests $F_1', F_2', \dots, F_m', F_{m+1}'$, again not having a cooperative coloring.

  The vertex set common to the new forests is $V' = (V \cup \sset{z}) \times V$, namely the vertex set consists of $\abs{V}+1$ copies of $V$. For every $u \in V\cup \sset{z}$ and every $i \in [m]$, take a copy $F^u_i$ of $F_i$ on the vertex set $\dset{(u,v)}{v \in V}$. Let $F_i'$ consist of $\abs{V}+1$ disjoint copies of $F_i$:
  $$F'_i := \bigcup_{u\in V\cup\sset{z}} F^u_i,\quad\text{ for all }i \in [m].$$
  To these we add the $(m+1)$st forest $F'_{m+1}$ obtained by joining $(z, u)$ to $(u, v)$ for all $u, v \in V$. So $F'_{m+1}$ is a disjoint union of stars, each with $\abs{V}$ leaves.

  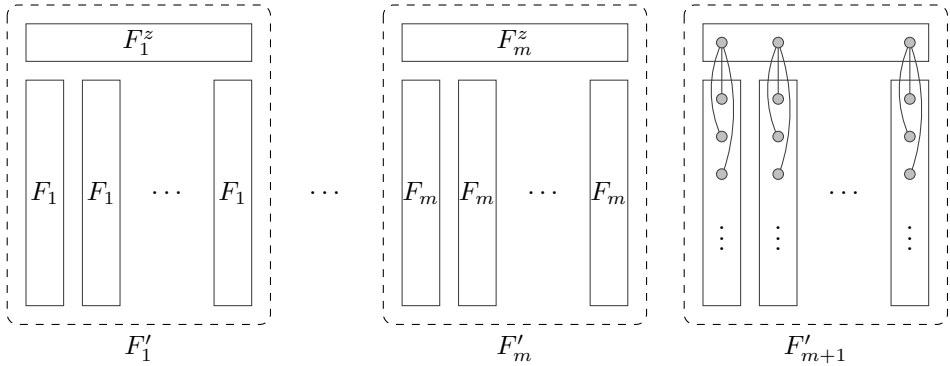
\begin{figure}[b]\centering
    \begin{tikzpicture}
      \node[block_h] at (0,2) {$F_1^z$};
      \node[block] at (-1.25,0) {$F_1$};
      \node[block] at (-0.5,0) {$F_1$};
      \node[block] at (1.25,0) {$F_1$};
      \path[draw, dashed, rounded corners] (-1.75,2.5) -- (-1.75,-1.75) -- node[below] {\small $F_1'$} (1.75,-1.75) -- (1.75,2.5) -- cycle;
      \node[block_h] at (5+0,2) {$F_m^z$};
      \node[block] at (5-1.25,0) {$F_m$};
      \node[block] at (5-0.5,0) {$F_m$};
      \node[block] at (5+1.25,0) {$F_m$};
      \path[draw, dashed, rounded corners] (5-1.75,2.5) -- (5-1.75,-1.75) -- node[below] {\small $F_m'$} (5+1.75,-1.75) -- (5+1.75,2.5) -- cycle;
      \node[block_h] at (9+0,2) {};
      \node[block] at (9-1.25,0) {};
      \node[block] at (9-0.5,0) {};
      \node[block] at (9+1.25,0) {};
      \path[draw, dashed, rounded corners] (9-1.75,2.5) -- (9-1.75,-1.75) -- node[below] {\small $F_{m+1}'$} (9+1.75,-1.75) -- (9+1.75,2.5) -- cycle;
      \node at (0.4,0) {\small$\dots$};
      \node at (2.5,0) {\small$\dots$};
      \node at (5.4,0) {\small$\dots$};
      \node at (9.4,0) {\small$\dots$};
      \node at (7.75,-0.5) {\small$\vdots$};
      \node at (8.5,-0.5) {\small$\vdots$};
      \node at (10.25,-0.5) {\small$\vdots$};

      \draw[darkgray] (7.75,2)--(7.75,1.25);
      \draw[darkgray] (7.75,2) arc (155:205:1.5);
      \draw[darkgray] (7.75,2) arc (20:-20:2.6);

      \draw[darkgray] (8.5,2)--(8.5,1.25);
      \draw[darkgray] (8.5,2) arc (155:205:1.5);
      \draw[darkgray] (8.5,2) arc (20:-20:2.6);

      \draw[darkgray] (10.25,2)--(10.25,1.25);
      \draw[darkgray] (10.25,2) arc (155:205:1.5);
      \draw[darkgray] (10.25,2) arc (20:-20:2.6);

      \draw[darkgray] (7.75,2) node[vertex]{};
      \draw[darkgray] (8.5,2) node[vertex]{};
      \draw[darkgray] (10.25,2) node[vertex]{};
      \draw[darkgray] (7.75,1.25) node[vertex]{};
      \draw[darkgray] (7.75,0.75) node[vertex]{};
      \draw[darkgray] (7.75,0.25) node[vertex]{};
      \draw[darkgray] (8.5,1.25) node[vertex]{};
      \draw[darkgray] (8.5,0.75) node[vertex]{};
      \draw[darkgray] (8.5,0.25) node[vertex]{};
      \draw[darkgray] (10.25,1.25) node[vertex]{};
      \draw[darkgray] (10.25,0.75) node[vertex]{};
      \draw[darkgray] (10.25,0.25) node[vertex]{};
    \end{tikzpicture}
    \caption{Construction of $Q(\cs) = (F_1', \dots, F_m', F_{m+1}')$ from $\cs = (F_1, \dots, F_m)$.} \label{fig2}
  \end{figure}
  Assume that there is a cooperative coloring $(I_1, I_2, \ldots, I_m, I_{m+1})$ for the system $Q(\cs)$. Since the forests $F_1^u, F_2^u, \dots ,F_m^u$ do not have a cooperative coloring, $I_{m+1}$ must contain a vertex from $\sset{u} \times V$ for all $u\in V\cup\sset{z}$. In particular, $I_{m+1}$ contains a vertex $(z,u) \in I_{k+1}$ for some $u \in V$ and a vertex $(u, v)$ for some $v\in V$. Since $(z,u)$ is connected in $F'_{m+1}$ to $(u,v)$, this is contrary to our assumption that $I_{m+1}$ is independent.

  Note that $\abs{V'} = \abs{V}^2+\abs{V} \le 2\abs{V}^2$. Note also that the maximum degree of $Q(\cs)$ is attained in $F'_{m+1}$, and it is equal to $\abs{V}$. Recursively define the system $\ct_m := Q(\ct_{m-1})$ consisting of $m$ forests for $m \ge 3$. Because the base $\ct_2$ has $4$ vertices, one can check inductively that $\abs{V(\ct_m)}$ is at most $2^{3\cdot2^{m-2}-1}$ using $\abs{V(\ct_m)} \le 2\abs{V(\ct_{m-1})}^2$. Thus the maximum degree of $\ct_m$ is at most $2^{3\cdot2^{m-3}-1} \le 2^{2^{m-1}}$.

  Given the maximum degree $d \ge 2$, choose $m := \ceil{\log_2\log_2d}$. By the choice of $m$, the maximum degree of $\ct_m$ is at most $2^{2^{m-1}} \le d$. By adding a few edges between the leaves in each forest of $\ct_m$, we can obtain a system of $m$ trees of maximum degree $d$ that does not have a cooperative coloring. This means $m_\ct(d) > m > \log_2\log_2d$.
\end{proof}

\begin{proof}[Proof of the upper bound on $m_{\ct}(d)$]
  Let $(T_1, T_2, \ldots ,T_m)$ be a system of trees of maximum degree $d$. We shall find a cooperative coloring by a random construction if $m \ge (1+o(1))\log_{4/3}d$.

  Choose arbitrarily for each tree $T_i$ a root so that we can specify the parent or a sibling of a vertex that is not the root of $T_i$. For each $T_i$, choose independently a random vertex set $S_i$, in which each vertex is included in $S_i$ independently with probability $1/2$. Set $$R_i := \dset{v\in S_i}{\text{the parent of }v\text{ is not in } S_i \text{, or } v\text{ is a root}}.$$ Since among any two adjacent vertices in $T_i$ one is the parent of the other, $R_i$ is independent in $T_i$.

  We shall show that with positive probability the sets $R_i$ form a cooperative coloring. For each vertex $v$, let $B_v$ be the event that $v \not\in \bigcup_{i=1}^mR_i$. If $v$ is the root of $T_i$, then $\pr{v\in R_i} = 1/2$; otherwise $\pr{v\in R_i} = 1/4$. In any case, $\pr{v\not\in R_i} \le 3/4$, and so $\pr{B_v}\le (3/4)^m$. Notice that $B_v$ is only dependent on the events $B_u$ for $u$ that is the parent, a sibling or a child of $v$ in some $T_i$. Since the degree of $v$ is at most $d$, it follows that $B_v$ is dependent on less than $2md$ other events. By the symmetric version of the Lov\'asz Local Lemma (see for example \cite[Chapter~5]{MR3524748}), if
  \begin{equation}\label{lll_tree}
    e \times \left(\frac{3}{4}\right)^m \times 2md \le 1,
  \end{equation}
  then with positive  probability no $B_v$ occurs, meaning that the sets $R_i$ form a cooperative coloring. The inequality \eqref{lll_tree} indeed holds under the assumption that $m \ge (1+o(1))\log_{4/3} d$.
\end{proof}

\section{Bipartite graphs} \label{sec:bipartite}

\begin{proof}[Proof of the lower bound on $m_{\cb}(d)$]
  Given $d$, take $m = \ceil{\log_2 d}$. Let the vertex set be $\sset{0,1}^m$, and for $j\in [m]$ let $G_j$ be the complete bipartite graph between $V_j^0$ and $V_j^1$ where \[
    V_{j}^k = \dset{v\in\sset{0,1}^m}{v_j = k},\quad\text{for }k\in\sset{0,1}.
  \]
  Note that the degree of $G_j$ is $2^{m-1} \le d$.

  Suppose that $I_1, \dots, I_m$ are independent sets in $G_1, \dots, G_m$ respectively. As each $G_j$ is a complete bipartite graph, $I_j \subseteq V_{j}^{k_j}$ for some $k_j\in\sset{0,1}$. Thus $(1-k_1, \dots, 1-k_m)$ is not in any $I_j$, and so $I_1, \dots, I_m$ do not form a cooperative coloring. This means $m_\cb(d) > m \ge \log_2d$.
\end{proof}

\begin{proof}[Proof of the upper bound on $m_{\cb}(d)$]
  Let $\cg=(G_1, \ldots ,G_m)$ be a system of bipartite graphs on the same vertex set $V$ with maximum degree $d$. By a semi-random construction, we shall find a cooperative coloring if $m \ge (1+\varepsilon)\frac{2d}{\ln d}$ for fixed $\varepsilon > 0$ and $d$ sufficiently large. We may assume that $m = O(d)$ because of \eqref{basic}.

  For each $j \in [m]$, let $(L_j,R_j)$ be a bipartition of $G_j$. Define $J_L(v) := \dset{j\in[m]}{v\in L_j}$ and $J_R(v) := \dset{j\in[m]}{v\in R_j}$ for each vertex $v\in V$, and let $A := \dset{v\in V}{\abs{J_L(v)} \ge m/2}$. Set $B := V \setminus A$. Clearly, we have
  \begin{subequations}
    \begin{align}
      \abs{J_L(a)} \ge m/2,&\quad\text{for all }a\in A; \label{jl} \\
      \abs{J_R(b)} \ge m/2,&\quad\text{for all }b\in B. \label{jr}
    \end{align}
  \end{subequations}
  Consider the following random process.
  \begin{enumerate}[nosep]
    \item For each $a \in A$, choose $j = j(a) \in J_L(a)$ uniformly at random, and put $a$ in the set $I_j$.
    \item For each $b \in B$, choose arbitrarily $j \in J_R(b)\setminus\dset{j(a)}{a\in A, (a,b)\in E(G_j)} =: J_R'(b)$ as long as it is possible, and put $b$ in the set $I_j$. \label{step2}
  \end{enumerate}
  For any $a, a' \in A \cap I_j$, $a, a' \in L_j$ and so $(a,a')\not\in G_j$. This means $A \cap I_j$ is independent, and similarly $B \cap I_j$ is independent. For any $b\in B\cap I_j$ and $(a,b)\in E(G_j)$, by the definition of $J_R'(b)$, $j(a) \neq j$ and so $a\not\in I_j$. Therefore $I_j$ is independent for all $j\in [m]$.

  To prove the existence of a cooperative coloring it suffices to show that $J_R'(b)$ is nonempty for all $b\in B$ with positive probability. For a vertex $b \in B$, let $E_b$ be the contrary event, that is, the event that $J_R'(b)$ is empty.

  For a fixed $b\in B$, let us estimate from above the probability of $E_b$. For every $j \in J_R(b)$, let $E^j$   be the event that $j\not\in J'_R(b)$, that is the event that $j(a) = j$ for some $a \in A$ that is a neighbor of $b$ in $G_j$. For each $a\in A$ that is a neighbor of $b$ in $G_j$, we have \[
    \pr{j(a) = j} = \frac{1}{\abs{J_L(a)}} \stackrel{\eqref{jl}}{\le} \frac{2}{m} \le \frac{\ln d}{(1+\varepsilon)d}.
  \]
  As there are at most $d$ neighbors of $b$ in $G_j$, we have for sufficiently large $d$ that
  \begin{equation}\label{ejb}
    1 - \pr{E^j} \ge \left(1-\frac{\ln d}{(1+\varepsilon)d}\right)^d \ge \exp\left(-(1-\varepsilon)\ln d\right) = d^{\varepsilon-1} \ge \frac{8\ln d}{m}.
  \end{equation}

  We claim that the events $E^j$, $j \in J_R(b)$, are negatively correlated. This is easier to see with the complementary events $\bar{E}^j$, $j \in J_R(b)$. We have to show that for any choice of indices $j_1, \ldots ,j_t \in J_R(b)$ there holds
  \[
    \cpr{E^j}{\bar{E}^{j_1} \cap \bar{E}^{j_2}\cap \ldots \cap \bar{E}^{j_t}} \ge \pr{E^j}.
  \]
  The event $\bar{E}^{j_1} \cap \bar{E}^{j_2} \cap \dots \cap \bar{E}^{j_t}$ means that for all $a\in A$ if $a$ is a neighbor of $b$ in $G_{j_i}$ then $j(a)\neq j_i$. Then, for any $j \not\in \sset{j_1, \dots ,j_t}$, for those vertices $a\in A$ that are neighbors of $b$ in $G_j$, knowing that $j(a) \neq j_i$ for certain $i\in[t]$ increases the probability that $j(a) = j$, and therefore increases the probability of ${E^j}$.

  By the claim, the inequality \eqref{ejb} and the fact that $E_b = \bigcap_{j \in J_R(b)} E^j$, we have
  \[
    \pr{E_b} \le \prod_{j\in J_R(b)}\pr{E^j} \stackrel{\eqref{jr}}{\le} \left(1-\frac{8\ln d}{m}\right)^{\frac{m}{2}} \le \exp\left(-\frac{8\ln d}{m}\cdot\frac{m}{2}\right) = \frac{1}{d^4}.
  \]

  The event $E_b$ is dependent on less than $md^2$ other events $E_{b'}$, since for such dependence to exist it is necessary that $b'\in B$ is at distance at most $2$ from $b$ in some graph $G_j$. Thus, by the Lov\'asz Local Lemma, for the positive probability that none of $E_b$ occurs it suffices that
  $$e \times \frac{1}{d^4}\times md^2 \le 1,$$
  which indeed holds for $d$ sufficiently large as $m = O(d)$.
\end{proof}

\section{Cooperative covers} \label{sec:covers}

Cooperative coloring of graphs is a special case of a more general concept.

\begin{definition}
  Given a system $(C_1, \ldots, C_n)$ of (abstract) simplicial complexes, all sharing the same vertex set $V$, a {\em cooperative cover} is a choice of faces $f_i \in C_i$ such that $\bigcup_{i = 1}^n f_i=V$.
\end{definition}

A cooperative coloring for $(G_1, \ldots ,G_n)$ is the special case in which $C_i$ is the independence complex $I(G_i)$ of $G_i$, that is, the collection of all independent sets in $G_i$.

\begin{definition}
  Given a hypergraph $C$ with vertex set $V$, the \emph{edge covering number} $\rho(C)$ is the minimal number of hyperedges from $C$ whose union is $V$. For a class $\cc$ of simplicial complexes, let $n_\cc(b)$ denote the minimal number $n$, such that every system $(C_1, \dots, C_n)$ of simplicial complexes belonging to $\cc$ on the same vertex set $V$ satisfying $\rho(C_i) \le b$ for all $i \le n$,
  has a cooperative cover. Let $n_\cc(b)=\infty$ if no such $n$ exists.
\end{definition}

For example, consider the class $\ci$ of all the independence complexes of graphs. If $G$ is bipartite, then $\rho(I(G))\le 2$. Hence the fact that $m_\cb(d) \ge \log_2(d)$ for all $d \ge 2$ (see Theorem \ref{main}) implies $n_\ci(2)=\infty$.

There are natural classes $\cc$ of hypergraphs for which $n_\cc$ is finite. One of these is the class of simplicial complexes associated to polymatroids, as introduced in \cite{MR0270945}. A \emph{polymatroid} $(V, r)$ is defined via a rank function $r\from 2^V \to \N$, that is submodular, monotone increasing and is $0$ on the empty set. A \emph{$k$-polymatroid} is a polymatroid in which every singleton set has rank at most $k$. For example, a $k$-uniform hypergraph $H$ endowed with the function $r(E)=\abs{\cup E}$, for every subset of hyperedges $E$ in $H$, is a $k$-polymatroid.

Following the notation in \cite[Section~11]{MR859549}, given a $k$-polymatroid $(V, r)$, a set $M \subseteq V$ is called a \emph{matching} if $r(M) = k\abs{M}$. By the submodularity of the rank function $r$,  the matchings in a $k$-polymatroid form a simplicial complex on $V$, which we call the \emph{matching complex} of a $k$-polymatroid.

\begin{theorem} \label{coop_cover}
  Let $\cm_k$ be the class of all the matching complexes of $k$-polymatroids. Then $n_{\cm_k}(b) \le kb$ for every $b$.
\end{theorem}

The proof uses the (homotopic) connectivity $\eta(C)$ of a complex $C$. We refer to \cite[Section~2]{MR2231877} for background. We shall use the following two topological tools. Given a complex $C$ on $V$ and $U \subseteq V$, we denote by $C[U]$ the simplicial subcomplex induced on $U$.

\begin{theorem}[Topological Hall's theorem]\label{topological_hall}
  Let $C$ be a simplicial complex on the vertex set $V$ and let $\bigcup_{i=1}^m W_i$ be a partition of $V$. If for all $I \subseteq [m]$
  $$\eta\left(C\left[\bigcup_{i \in I}W_i\right]\right) \ge \abs{I},$$
  then $C$ contains a face $\sigma$ such that $\abs{\sigma \cap W_i} = 1$ for all $i \in [m]$.
\end{theorem}

\begin{theorem}\label{ab_unpublished}
  If $C$ is a matching complex on $V$ of a $k$-polymatroid, then the connectivity $\eta(C)$ of $C$ is at least $\nu(C)/k$, where  $\nu(C)$ is the maximal size of faces in $C$.
\end{theorem}

The above formulation of Theorem~\ref{topological_hall} first appeared in \cite{MR1805715}, attributed to the first author of the present paper (see the remark after Theorem 1.3 in \cite{MR1805715}). Theorem~\ref{ab_unpublished} is an unpublished result of the first two authors. The special case, where the $k$-polymatroid is the sum of $k$ matroids on the same vertex set, is proved in \cite[Theorem 6.5]{MR2231877}.

\begin{proof}[Proof of Theorem~\ref{coop_cover}]
  Let $n=kb$, and let $C_1, \dots, C_n$ be simplicial complexes associated to $k$-polymatroids $(V, r_1), \dots, (V, r_n)$ on the same vertex set $V$ such that the edge covering number of each $C_i$ is at most $b$.
  Let $C$ be the join of $C_1, \dots, C_n$ on $V\times[n]$, that is, $$C := \dset{\bigcup_{i=1}^n \sigma_i\times \sset{i}}{\sigma_i \in C_i \text{ for all }i\in[n]}.$$ A cooperative cover can be viewed as a face $\sigma \in C$ such that $\abs{\sigma \cap \left(\sset{v}\times [n]\right)} = 1$ for all $v \in V$. By the topological Hall's theorem, it suffices to prove that $$\eta\left(C[U\times [n]]\right) \ge \abs{U} \text{ for all }U \subseteq V.$$

  Let $U$ be a subset of $V$. Note that $C_i[U]$ is the matching complex of the $k$-polymatroid $(U, r_i|_{U})$. By Theorem~\ref{ab_unpublished}, $\eta(C_i[U]) \ge \nu(C_i[U]) / k$. Since $\nu(C_i[U])$ is the maximal size of faces in $C_i[U]$ and the edge covering number of $C_i[U]$ is at most $b$, we obtain $\nu(C_i[U])b \ge \abs{U}$, and so $\eta(C_i[U]) \ge \abs{U}/(kb)$.
  Notice that $C[U \times [n]]$ is the join of $C_1[U], \dots, C_n[U]$. Using the superadditivity of $\eta$ with respect to the join operator and Theorem~\ref{ab_unpublished}, we obtain the required condition for the topological Hall's theorem
  \begin{equation*}
    \eta(C[U\times [n]]) \ge \sum_{i=1}^n \eta(C_i[U]) \ge \sum_{i=1}^n {\abs{U}}/(kb) = \abs{U}.\qedhere
  \end{equation*}
\end{proof}

\begin{remark}
  It is of interest to explore the sharpness of this result.
\end{remark}

\bibliographystyle{alpha}
\bibliography{cooperative_coloring}

\begin{thebibliography}{AHHS15}

\bibitem[AB06]{MR2231877}
Ron Aharoni and Eli Berger.
\newblock The intersection of a matroid and a simplicial complex.
\newblock {\em Trans. Amer. Math. Soc.}, 358(11):4895--4917, 2006.

\bibitem[ABZ07]{MR2345810}
Ron Aharoni, Eli Berger, and Ran Ziv.
\newblock Independent systems of representatives in weighted graphs.
\newblock {\em Combinatorica}, 27(3):253--267, 2007.

\bibitem[AHHS15]{MR3359930}
Ron Aharoni, Ron Holzman, David Howard, and Philipp Spr\"ussel.
\newblock Cooperative colorings and independent systems of representatives.
\newblock {\em Electron. J. Combin.}, 22(2):Paper 2.27, 14, 2015.

\bibitem[AK98]{MR1774970}
Noga Alon and Michael Krivelevich.
\newblock The choice number of random bipartite graphs.
\newblock {\em Ann. Comb.}, 2(4):291--297, 1998.

\bibitem[AR08]{MR2570629}
Omid Amini and Bruce Reed.
\newblock List colouring constants of triangle free graphs.
\newblock In {\em The {IV} {L}atin-{A}merican {A}lgorithms, {G}raphs, and
  {O}ptimization {S}ymposium}, volume~30 of {\em Electron. Notes Discrete
  Math.}, pages 135--140. Elsevier Sci. B. V., Amsterdam, 2008.

\bibitem[AS16]{MR3524748}
Noga Alon and Joel~H. Spencer.
\newblock {\em The probabilistic method}.
\newblock Wiley Series in Discrete Mathematics and Optimization. John Wiley \&
  Sons, Inc., Hoboken, NJ, fourth edition, 2016.

\bibitem[Edm70]{MR0270945}
Jack Edmonds.
\newblock Submodular functions, matroids, and certain polyhedra.
\newblock In {\em Combinatorial {S}tructures and their {A}pplications ({P}roc.
  {C}algary {I}nternat. {C}onf., {C}algary, {A}lta., 1969)}, pages 69--87.
  Gordon and Breach, New York, 1970.

\bibitem[Hax01]{MR1860440}
P.~E. Haxell.
\newblock A note on vertex list colouring.
\newblock {\em Combin. Probab. Comput.}, 10(4):345--347, 2001.

\bibitem[LP86]{MR859549}
L.~Lov\'{a}sz and M.~D. Plummer.
\newblock {\em Matching theory}, volume 121 of {\em North-Holland Mathematics
  Studies}.
\newblock North-Holland Publishing Co., Amsterdam; North-Holland Publishing
  Co., Amsterdam, 1986.
\newblock Annals of Discrete Mathematics, 29.

\bibitem[LS07]{MR2354708}
Po-Shen Loh and Benny Sudakov.
\newblock Independent transversals in locally sparse graphs.
\newblock {\em J. Combin. Theory Ser. B}, 97(6):904--918, 2007.
\newblock \href{https://arxiv.org/abs/0706.2124}{\tt arXiv:0706.2124[math.CO]}.

\bibitem[Mes01]{MR1805715}
Roy Meshulam.
\newblock The clique complex and hypergraph matching.
\newblock {\em Combinatorica}, 21(1):89--94, 2001.

\end{thebibliography}

\end{document}